# FINITE-DIMENSIONAL REALIZATIONS FOR STOCHASTIC EQUATIONS IN THE HJM-FRAMEWORK

DAMIR FILIPOVIĆ AND JOSEF TEICHMANN

ABSTRACT. This paper discusses finite-dimensional (Markovian) realizations (FDRs) for Heath–Jarrow–Morton interest rate models. We consider a $d$-dimensional driving Brownian motion and stochastic volatility structures that are non-degenerate smooth functionals of the current forward rate. In a recent paper, Björk and Svensson give sufficient and necessary conditions for the existence of FDRs within a particular Hilbert space setup. We extend their framework, provide new results on the geometry of the implied FDRs and classify all of them. In particular, we prove their conjecture that every short rate realization is 2-dimensional. More generally, we show that all generic FDRs are at least $(d+1)$-dimensional and that all generic FDRs are affine. As an illustration we sketch an interest rate model, which goes well with the Svensson curve-fitting method. These results cannot be obtained in the Björk–Svensson setting.

A substantial part of this paper is devoted to analysis on Fréchet spaces, where we derive a Frobenius theorem. Though we only consider stochastic equations in the HJM-framework, many of the results carry over to a more general setup.

## 1. INTRODUCTION

This paper discusses the existence of finite-dimensional forward rate models within the framework of Heath, Jarrow and Morton (henceforth HJM) [10]. Following Musiela [15], we denote by $r_t(x)$ the instantaneous forward rate at time $t$ for maturity $t + x$. The price at time $t$ of a zero coupon bond with maturity $T$ is then given by

$$P(t, T) = \exp\left(-\int_0^{T-t} r_t(x)\, dx\right).$$

It is shown in [6] that essentially any HJM-type model for $r_t(x)$ can be considered as a stochastic equation

$$\begin{cases} dr_t = (Ar_t + \alpha_{HJM}(t, r_t))\, dt + \sigma(t, r_t)\, dW_t \\ r_0 = h_0 \end{cases} \tag{1.1}$$

on a separable Hilbert space $H$ which is characterized by the following three properties:

(H1) $H \subset C(\mathbb{R}_{\geq 0}; \mathbb{R})$ with continuous embedding[1]; that is, the pointwise evaluations $ev_x : h \mapsto h(x)$ are continuous linear functionals.

---

[1] Typically, $H$ consists of equivalence classes of functions, which we shall identify with their continuous representatives, respectively.

(H2) The family of right-shifts, $S(t)f(x) = f(t+x)$, forms a strongly continuous semigroup $\{S(t) \mid t \in \mathbb{R}_{\geq 0}\}$ on $H$.

(H3) There exists a closed subspace $H_0$ of $H$ such that

$$\mathcal{S}(f,g)(x) := f(x) \int_0^x g(\eta) \, d\eta,$$

defines a continuous bilinear mapping $\mathcal{S} : H_0 \times H_0 \to H$.

We write shortly $\mathcal{S}(f)$ for $\mathcal{S}(f,f)$.

As an illustration we shall always have the following example in mind (see [6, Section 5]).

**Example 1.** *Let $w : \mathbb{R}_{\geq 0} \to [1, \infty)$ be a non-decreasing $C^1$-function such that $w^{-1/3} \in L^1(\mathbb{R}_{\geq 0})$. We may think of $w(x) = e^{\alpha x}$ or $w(x) = (1+x)^\alpha$, for $\alpha > 0$ or $\alpha > 3$, respectively. The space $H_w$ consisting of absolutely continuous functions $h$ on $\mathbb{R}_{\geq 0}$ and equipped with the norm*

$$\|h\|_w^2 := |h(0)|^2 + \int_{\mathbb{R}_{\geq 0}} |h'(x)|^2 w(x) \, dx$$

*is a Hilbert space which satisfies (H1)–(H2). Property (H3) is satisfied for $H_0 = H_w^0 := \{h \in H_w \mid \lim_{x \to \infty} h(x) = 0\}$.*

Throughout we are given a filtered probability space $(\Omega, \mathcal{F}, (\mathcal{F}_t), \mathbb{P})$ satisfying the usual conditions. The predictable $\sigma$-field is denoted by $\mathcal{P}$. Now we can give a meaning to the coefficients in (1.1). Here $W = (W^1, \ldots, W^d)$, $d \in \mathbb{N}$, denotes a standard $d$-dimensional Brownian motion. The volatility coefficient $\sigma(t, \omega, h) = (\sigma^1(t, \omega, h), \ldots, \sigma^d(t, \omega, h))$ is a measurable mapping from $(\mathbb{R}_{\geq 0} \times \Omega \times H, \mathcal{P} \otimes \mathcal{B}(H))$ into $(H_0^d, \mathcal{B}(H_0^d))$, well-defining

$$\alpha_{HJM}(t, \omega, h) := \sum_{j=1}^d \mathcal{S}(\sigma^j(t, \omega, h)),$$

a measurable mapping from $(\mathbb{R}_{\geq 0} \times \Omega \times H, \mathcal{P} \otimes \mathcal{B}(H))$ into $(H, \mathcal{B}(H))$.

The operator $A$ is the infinitesimal generator of the semigroup $\{S(t) \mid t \in \mathbb{R}_{\geq 0}\}$. Its domain in $H$ is denoted by $D(A)$, and we write $A^*$ for its adjoint. Both operators, $A$ and $A^*$, are closed. It is easy to see that $D(A) \subset \{h \in H \cap C^1(\mathbb{R}_{\geq 0}; \mathbb{R}) \mid h' \in H\}$ and $Ah = h'$. Without much loss of generality we shall in fact assume

(H4) $D(A) = \{h \in H \cap C^1(\mathbb{R}_{\geq 0}; \mathbb{R}) \mid h' \in H\}$.

Also (H4) is satisfied for the spaces $H_w$ from Example 1.

We distinguish, in decreasing order of generality, between (local) mild, weak and strong solutions to equation (1.1) and its time-shifted versions. The reader is referred to [5] or [6] for the precise definitions.

**Definition 1.** *A subset $V$ of $H$ is called* locally invariant *for (1.1) if, for all space-time initial points $(t_0, h_0) \in \mathbb{R}_{\geq 0} \times V$, there exists a continuous local weak solution $r^{(t_0, h_0)}$ to (1.1) with lifetime $\tau = \tau(t_0, h_0)$ such that*

$$r_{t \wedge \tau}^{(t_0, h_0)} \in V, \quad \forall t \in \mathbb{R}_{\geq 0}.$$

For the notion of a finite-dimensional $C^k$-submanifold of a Hilbert space we refer to [13]. We restate [8, Theorem 3].

**Theorem 1.** *Suppose that $\sigma(t, \omega, h)$ is continuously differentiable in $h$ and right-continuous in $t$. Let $\mathcal{M}$ be an $m$-dimensional $C^2$-submanifold of $H$. Then the following conditions are equivalent:*

i) $\mathcal{M}$ *is locally invariant for (1.1)*

ii) $\mathcal{M} \subset D(A)$ *and*

$$\mu(t, \omega, h) := Ah + \alpha_{HJM}(t, \omega, h) - \frac{1}{2} \sum_{j=1}^{d} D\sigma^j(t, \omega, h)\sigma^j(t, \omega, h) \in T_h\mathcal{M} \quad (1.2)$$

$$\sigma^j(t, \omega, h) \in T_h\mathcal{M}, \quad j = 1, \ldots, d, \quad (1.3)$$

*for all $(t, h) \in \mathbb{R}_{\geq 0} \times \mathcal{M}$, for $\mathbb{P}$-a.e. $\omega$.*

Hence the stochastic invariance problem to (1.1) is equivalent to the deterministic invariance problems (for any $\omega$ apart from a $\mathbb{P}$-nullset) related to the $t$-dependent vector fields $\{\mu, \sigma^1, \ldots, \sigma^d\}$.

Note that Theorem 1 provides conditions for the invariance of a single given submanifold $\mathcal{M}$. It does not say anything about the *existence* of finite-dimensional invariant submanifolds for (1.1) or $\{\mu, \sigma^1, \ldots, \sigma^d\}$, respectively. This issue will be exploited in the present paper.

We remark that the existence of finite-dimensional invariant submanifolds is essentially equivalent to the existence of finite-dimensional realizations *(FDR)* for (1.1) in the following sense (see [6, Theorem 6.4.1]).

**Theorem 2.** *Let $\sigma$ and $\mathcal{M}$ be as in Theorem 1. Suppose $\mathcal{M}$ is locally invariant for (1.1). Then, for any space-time initial point $(t_0, h_0) \in \mathbb{R}_{\geq 0} \times \mathcal{M}$, there exists a $C^2$-map $\phi : \mathbb{R}^m \to U \cap \mathcal{M}$, where $U$ is an open neighborhood of $h_0$ in $H$, and an $\mathbb{R}^m$-valued diffusion process $Z$ such that*

$$r_{t \wedge \tau}^{(t_0, h_0)} = \phi(Z_{t \wedge \tau}), \quad \forall t \in \mathbb{R}_{\geq 0}, \quad (1.4)$$

*for some stopping time $\tau > 0$.*

Following [3] we say that equation (1.1) has a local *$m$-dimensional realization* at $(t_0, h_0)$ if (1.4) holds. We say that (1.1) *generically* admits an FDR at $(t_0, h_0)$ if (1.1) has a local $m$-dimensional realization for initial points in a neighborhood of $(t_0, h_0)$.

**Remark 1.** *In view of the preceding remarks we shall suppress the $\omega$-dependence from now on, making the convention that all subsequent statements (have to) hold simultaneously for all $\omega$ (apart from a $\mathbb{P}$-nullset).*

**Remark 2.** *The forward rate is a mathematical idealization rather than an object that can directly be observed on the market. For statistical inference one usually estimates the current forward curve by fitting a parametrized family $\mathcal{G}$ of curves to the finitely given many market data. Such a family $\mathcal{G}$ is generically a finite-dimensional immersed submanifold in $H$. For the sake of consistency of the curve-fitting procedure with any given stochastic interest rate model of the form (1.1), it is vital that $\mathcal{G}$ is invariant for the flow of (1.1). If this is not the case, it is unclear what the parameters of the daily fitting should mean in terms of the model. Consequently, the choice of the volatility structure, $\sigma^1, \ldots, \sigma^d$, determines the consistent curve-fitting procedures. Volatility structures that do not admit invariant submanifolds (i.e. FDRs) should therefore be avoided.*

The FDR-problem consists of finding sufficient conditions on $\sigma^1, \ldots, \sigma^d$ for the existence of generic FDRs. Björk et al [1], [3] translated this into an appropriate geometric language. In [3] they completely solved the FDR-problem for a very particular choice of $H$. Their key argument is the classical Frobenius theorem, since they are looking for foliations (which is the appropriate notion for the FDR-problem on Hilbert spaces). Therefore they define a Hilbert space, $\mathcal{H}$, on which $A = d/dx$ is a bounded linear operator. As a consequence $\mathcal{H}$ consists solely of entire analytic functions (see [3, Proposition 4.2]). It is well known however that the forward rates implied by a Cox–Ingersoll–Ross (CIR) [4] short rate model are of the form $r_t(x) = g_0(x) + r_t(0)g_1(x)$ where

$$g_0(x) = d\,\frac{e^{ax} - 1}{e^{ax} + c} \quad \text{and} \quad g_1(x) = \frac{be^{ax}}{(e^{ax} + c)^2},$$

for some $a, b, c > 0$ and $d \geq 0$ (see e.g. [6, Section 7.4.1]). Since both $g_0$ and $g_1$, when extended to $\mathbb{C}$, have a singularity at $x = (\log(c) + i\pi)/a$, they cannot be entire analytic. Hence the CIR forward rates do not belong to $\mathcal{H}$! This fact seems to be inconsistent with some of the results in [3] (see Remark 7.4 therein). Since CIR-type structures appear generically in HJM-framework, the Björk-Svensson setting is too narrow, even though all geometric ideas are already formulated.

To overcome this difficulty we have to choose a larger forward curve space. However, we cannot get around the Frobenius theorem, which requires smoothness of the vector fields $\{\mu, \sigma^1, \ldots, \sigma^d\}$. The problem here is that $A$ is typically an unbounded operator (the choice of $\mathcal{H}$ in [3] is exactly made to overcome this problem). The appropriate framework for an extended version of the Frobenius theorem is given by the Fréchet space

$$D(A^\infty) := \bigcap_{n \in \mathbb{N}} D(A^n),$$

equipped with the family of seminorms

$$p_n(h) = \sum_{m=0}^{n} \|A^m h\|_H, \quad n \in \mathbb{N}_0.$$

Indeed, $D(A^\infty)$ is essentially the largest subspace of $H$ which is invariant for $A$. This calls for a Frobenius theorem in Fréchet spaces, which is definitively not straightforward, since the inverse function theorem fails in Fréchet spaces in general (see [12] for example). Denote by $A_0 : D(A_0) \to H_0$ the restriction of $A$ to $H_0$. That is, $D(A_0) = \{h \in D(A) \cap H_0 \mid Ah \in H_0\}$. The definition of the Fréchet space $D(A_0^\infty) = \cap_{n \in \mathbb{N}} D(A_0^n)$ is obvious. The next result follows immediately from (H1), (H3) and (H4).

**Lemma 1.** *For any $f, g \in D(A_0)$ we have $\mathcal{S}(f, g) \in D(A)$ and*

$$A\mathcal{S}(f, g) = \mathcal{S}(Af, g) + \mathcal{S}(f, Ag) + f\,ev_0(g).$$

*Hence $\mathcal{S} : D(A_0^\infty) \times D(A_0^\infty) \to D(A^\infty)$ is a continuous bilinear mapping.*

The preceding specifications for $\sigma$ are still too general for concrete implementations. We actually have the idea of $\sigma$ being sensitive with respect to functionals of the forward curve. That is, $\sigma^j(t, h) = \phi^j(t, \ell_1(h), \ldots, \ell_q(h))$, for some $q \in \mathbb{N}$, where $\phi^j : \mathbb{R}_{\geq 0} \times \mathbb{R}^q \to D(A_0^\infty)$ is a smooth map and $\ell_1, \ldots, \ell_q$ denote continuous linear functionals on $H$ (or even on $C(\mathbb{R}_{\geq 0}; \mathbb{R})$). We may think of $\ell_i(h) = (1/x_i) \int_0^{x_i} h(\eta) \, d\eta$ (benchmark yields) or $\ell_i(h) = ev_{x_i}(h)$ (benchmark forward rates). This idea is (generalized and) expressed in terms of the following regularity and non-degeneracy assumptions:

(A1)  $\sigma(t,h) = \phi(t, \ell(h))$ with $\ell : H \to B$ smooth linear and $\phi : B \to D(A_0^\infty)^d$ smooth, so it is a Banach-map (see Definition 3 below).

(A2)  $(\ell, \ell \circ A) : D(A^\infty) \to B \times B$ is a linear open map.

(A3)  $\{\phi^1(t,y), \ldots, \phi^d(t,y)\}$ are linearly independent, for all $(t,y) \in \mathbb{R}_{\geq 0} \times B$.

(A4)  $A$ is unbounded, that is, $D(A)$ is a strict subset of $H$.

We believe that this setup is flexible enough to capture any reasonable HJM-type interest rate model. Combining Lemma 1 and (A1) yields

**Lemma 2.** $\mathcal{S}(\sigma^j(\cdot)) : \mathbb{R}_{\geq 0} \times H \to D(A^\infty)$ is a Banach-map, for any $j = 1, \ldots, d$, hence also $\alpha_{HJM}$.

Assumption (A4) is equivalent to the following property of $A$.

**Lemma 3.** Let $A$ be the generator of a strongly continuous semigroup $S$ on a Banach space $Y$, then the operator $A : D(A^\infty) \to D(A^\infty)$ is a Banach-map if and only if $A : Y \to Y$ is bounded.

*Proof.* Let $Y$ be a Banach space and $A$ the generator of a strongly continuous semigroup $S$. If $A : D(A^\infty) \to D(A^\infty)$ is in a neighborhood of $U$ a point $x_0$ a Banach map, then there are smooth mappings $R : U \subset D(A^\infty) \to X$ and $Q : V \subset X \to D(A^\infty)$ such that $A = Q \circ R$ and $X$ is a Banach space. By differentiation we obtain at $x_0$

$$A = DQ(x_0) \cdot DR(x_0)$$

which means in particular by continuity that there exists $\underline{n \geq 0}$ such that $DR(x_0)$ can be extended continuously to a linear mapping $\overline{DR(x_0)} : D(A^n) \to X$. So $A : D(A^n) \to D(A^n)$ is a continuous mapping. We recall the Sobolev Hierarchy for strongly continuous semigroups (see [16]) defined by the following commutative diagram

Here $R(\lambda) := (\lambda - A)^{-1}$ denotes the resolvent at a point of the resolvent set, which defines an isomorphism from $D(A^n)$ to $D(A^{n+1})$. The semigroups $S^{(n)}$ are defined by restriction and are strongly continuous in the respective topologies. The generator of $S^{(n)}$ is $A$ restricted to $D(A^{(n+1)})$. If $A$ were continuous on $D(A^{(n+1)})$, then we could climb up the diagram by the resolvent isomorphisms. $A$ continuous on $D(A^{(n)})$ implies that $S^{(n)}$ is a smooth group, so by climbing up through the isomorphisms $S^{(0)}$ is a smooth group and therefore the infinitesimal generator is continuous, since it is everywhere defined, by the closed graph theorem. ◻

We will demonstrate that the Frobenius integrability condition implies a very particular geometry of the invariant submanifolds – loosely speaking, each of them is a band of copies of an affine submanifold. This is similar to Hamiltonian mechanics, where the existence of several conservation laws forces the flow to be viable not only on a low dimensional submanifold, but even on a special one, namely a torus.

**Remark 3.** *Although we focus on the particular stochastic equation (1.1), many arguments can be carried over to more general stochastic equations in the spirit of Da Prato and Zabczyk [5] (see also Remark 1).*

## 2. Analysis on Fréchet Spaces

A Fréchet space $E$ is a completely metrizable locally convex vector space. On a Fréchet space $E$ we are equivalently given a sequence (which can be chosen increasing) of seminorms $\{p_n\}_{n \geq 1}$ generating a complete locally convex topology. For the construction of a differential calculus we need the concept of differentiable curves into the Fréchet space and differentiable mappings on open subsets of a Fréchet space. To obtain a completely working differential calculus we should enlarge the category of vector spaces to so called convenient vector spaces, however, for our purposes we can directly give the definitions on Fréchet spaces (see [12] for all details). We remark that already on Fréchet spaces the situation concerning analysis was complicated and unclear until convenient calculus was invented (see [12] for extensive historical remarks).

We denote the set of continuous linear functionals by $E'$. A subset $B$ of a Fréchet space is bounded if and only if $l(B)$ is a bounded subset of $\mathbb{R}$ for all $l \in E'$. A multilinear mapping $m : E_1 \times ... \times E_n \to F$ is called bounded if bounded sets $B_1 \times ... \times B_n$ are mapped to bounded subsets of the Fréchet space $F$. On Fréchet spaces bounded linear mappings are continuous (see [12] or any textbook on locally convex spaces for the proof). In general this is not the case, however continuous linear functionals are clearly bounded. The set of bounded linear operators on a Fréchet space $E$ is denoted by $L(E)$. We prefer the notion "bounded" to "continuous" in the case of linear operators, since it is better adapted to convenient calculus (the following result is true for all convenient spaces).

**Definition 2.** *Let $E$ be a Fréchet space, then $c : \mathbb{R} \to E$ is called smooth if all derivatives exist as limits of difference quotients. The set of smooth curves is denoted by $C^\infty(\mathbb{R}, E)$. A mapping $f : U \subset E \to \mathbb{R}$ is called smooth if $U$ is open and $f \circ c \in C^\infty(\mathbb{R}, \mathbb{R})$ for all $c \in C^\infty(\mathbb{R}, E)$.*

**Theorem 3.** *Let $E, G, H$ be Fréchet spaces, $U \subset E$, $V \subset G$ open subsets, then we obtain:*

i) *A curve $c : \mathbb{R} \to E$ is smooth if and only if it is weakly smooth, i.e. $l \circ c \in C^\infty(\mathbb{R}, \mathbb{R})$ for all $l \in E'$.*

ii) *Multilinear mappings are smooth if and only if they are bounded.*

iii) *If $f : U \to G$ is smooth, then $df : U \times E \to G$ is smooth and bounded linear in the second component, where*

$$df(x, v) := \frac{d}{dt}|_{t=0} f(x + tv)$$

iv) *The chain rule holds.*
v) *Taylor's formula is true, where higher derivatives are defined as usual:*

$$f(x+y) = \sum_{i=0}^{n} \frac{1}{i!} d^i f(x) y^{(i)} + \int_0^1 \frac{(1-t)^n}{n!} d^{n+1} f(x+ty) \, (y^{(n+1)}) dt$$

*for all $n \in \mathbb{N}$.*
vi) *Given $c : \mathbb{R} \to L(E)$, then $c : \mathbb{R} \times E \to E$ is smooth if and only if $ev_v \circ c \in C^\infty(\mathbb{R}, E)$ for all $v \in E$, where $ev_v : L(E) \to E$ denotes the point evaluation.*
vii) *A family $a_i \in L(E)$ for $i \in I$ of bounded linear operators is uniformly bounded, i.e. for all bounded subsets $B \subset E$ there is a bounded subset $C \subset E$ such that for all $i \in I$ we have $a_i(B) \subset C$, if and only if it is pointwise bounded ("smooth uniform boundedness principle").*
viii) *Let $B$ be a Banach space and $f : U \to G$ be a smooth mapping, then $f$ is locally Lipschitz, i.e. for every point $x \in U$ there is a neighborhood $V \subset U$ with*

$$\{\frac{f(x) - f(y)}{\|x - y\|} | x \neq y \in V\} \text{ is bounded in } G$$

**Remark 4.** *The introduction of analysis on Fréchet spaces via definition 2 is not limited to Fréchet spaces, but works on all locally convex spaces (so called convenient vector spaces) where weakly smooth curves are smooth. The necessity to work on these spaces instead of working on Fréchet spaces arises from the fact that spaces like $C^\infty(E, F)$ or even $L(E)$ generically are not Fréchet space for Fréchet spaces $E, F$.*

*Even in the case of $\mathbb{R}^2$ it is not obvious that the smooth functions of definition 2 are smooth in the classical sense. Notice that convenient calculus completely solves the problem, how to do analysis on locally convex vector spaces, and provides a very powerful tool in concrete calculations (see [12] for many examples).*

Notice that on Fréchet spaces the notion of $C^n$-curves is useful, too, and that we are also given the chain rule with smooth functions by the uniform boundedness principle:

$$\frac{f(c(t)) - f(c(0))}{t} = \int_0^1 Df(c(0) + s(c(t) - c(0))) \cdot \frac{c(t) - c(0)}{t} ds$$
$$\longrightarrow Df(c(0)) \cdot c'(0)$$

as $t \to 0$. Furthermore we obtain the following useful result by Taylor's formula and uniform boundedness (see [9] for details).

**Theorem 4.** *Let $E, G$ be Fréchet spaces and $f : U \subset E \to G$ a smooth mapping. Given $x_0 \in U$ and a continuous seminorm $q$ on $G$, then there is a continuous seminorm $p$ on $E$ and $\varepsilon > 0$ such that*

$$q(f(x) - f(y)) \leq p(x - y)$$

*if $p(x - x_0) < \varepsilon$ and $p(y - x_0) < \varepsilon$.*

Concerning differential equations, there are possible counterexamples on non-normable Fréchet spaces in all directions, which causes some problems in the foundations of differential geometry (see [12] and the excellent review article [14]). Nevertheless a useful generalization of the existence theorem for differential equations on Banach spaces is given by the following Banach mapping principle (see [9] for details). Given $P : U \subset E \to E$ smooth, we are looking for solutions of the ordinary differential equation

$$x :] - \varepsilon, \varepsilon[ \to U \text{ smooth}$$

$$\frac{d}{dt} x(t) = P(x(t))$$

$$x(0) \in U$$

If for any initial value $y$ in a small neighborhood $V$ of $x_0 \in U$ there is a unique smooth solution $t \mapsto x_y(t)$ for $t \in ] - \varepsilon, \varepsilon[$ depending smoothly on $y$, then $Fl(t, y) := x_y(t)$ defines a local flow, i.e. a smooth map

$$Fl :] - \varepsilon, \varepsilon[ \times V \to E$$

$$Fl(0, y) = y$$

$$Fl(t, Fl(s, y)) = Fl(s + t, y)$$

if $s, t, s + t \in ] - \varepsilon, \varepsilon[$ and $Fl(s, y) \in V$. If there is a local flow around $x_0 \in U$ (this shall mean once and for all: "in an open neighborhood of $x_0$"), then the differential equation is uniquely solvable around $x_0 \in U$ and the dependence on initial values is smooth. Let $U$ be connected, then a maximal local flow $Fl$ associated to $P$ is a local flow defined on a connected open set $W \subset \mathbb{R} \times U$ such that for any integral curve $x : I \to U$ with $x(0) = y$ and $I$ connected we have $I \times \{y\} \in W$ and $Fl(t, y) = x(t)$ for $t \in I$. If $P$ admits a local flow around any point $x_0 \in U$, then it admits a maximal local flow. A maximal local flow satisfies

$$(t, Fl(s, y)) \in W \text{ if and only if } (s + t, y) \in W$$

for all $(s, y) \in W$ and $t \in \mathbb{R}$. In this case $Fl(t, Fl(s, y)) = Fl(s + t, y)$ (on flows see for example [12]).

**Definition 3.** *Given a Fréchet space $E$, a smooth mapping $P : U \subset E \to E$ is called a Banach map if there are smooth (not necessarily linear!) mappings $R : U \subset E \to X$ and $Q : V \subset X \to E$ such that $P = Q \circ R$.*

$$U \subset E \xrightarrow{\quad P \quad} E$$

with $R$ downward to $V \subset X$ and $Q$ upward to $E$.

**Theorem 5.** *We denote by $\mathcal{B}(U)$ the set of Banach map vector fields, by $\mathfrak{X}(U)$ all vector fields on an open subset of the Fréchet space $E$. Then $\mathcal{B}(U)$ is a $C^\infty(U, \mathbb{R})$-submodule of $\mathfrak{X}(U)$.*

*Proof.* We have to show that for $f, g \in C^\infty(U, \mathbb{R})$ and $P_1, P_2 \in \mathcal{B}(U)$ the linear combination $fP_1 + gP_2 \in \mathcal{B}(U)$. Given $P_i = Q_i \circ R_i$ for $i = 1, 2$ with intermediate Banach space $X_i$, then $fP_1 + gP_2 = Q \circ R$ with $Q : \mathbb{R}^2 \times V_1 \times V_2 \subset \mathbb{R}^2 \times X_1 \times X_2 \to E$ and $R : U \to \mathbb{R}^2 \times X_1 \times X_2$ such that

$$Q(r, s, v_1, v_2) = r Q_1(v_1) + s Q_2(v_2)$$

$$R(x) = (f(x), g(x), R_1(x), R_2(x))$$

So the sum $fP_1 + gP_2$ is a Banach map and therefore the set of all Banach map vector fields carries the asserted submodule structure.                    □

**Theorem 6** (Banach mapping principle). *Let $P : U \subset E \to E$ be a Banach map, then $P$ admits a maximal local flow.*

*Proof.* For the proof see [9]. $\qquad\qquad\qquad\qquad\qquad\qquad\qquad\qquad\qquad\qquad\qquad$ $\square$

**Remark 5.** *Parameters and time-dependence are treated in the following way. Given an open subset of parameters $Z \subset Y$ of a Banach space and $P : I \times Z \times U \to E$, where $I$ is a open set in $\mathbb{R}$ and $U$ is open in $E$, such that $P_{t,p} = Q_{t,p} \circ R_{t,p}$, where $Q$ and $R$ depend smoothly on time and parameters, $P$ admits a unique smooth solution for any initial value $x_0 \in U$ at any time point $t_0 \in I$ depending smoothly on parameters, time and initial values. For the proof we look at the extended space $G := \mathbb{R} \times Y \times E$ with $\widetilde{P}(t, p, x) = (1, 0, P(t, p, x))$ and*

$$\widetilde{Q}(t, p, z) = (1, 0, Q_{t,p}(x))$$
$$\widetilde{R}(t, p, x) = (t, p, R_{t,p}(x))$$

*with Banach space $\widetilde{X} := \mathbb{R} \times Y \times X$.*

We can replace in the above definition the interval $]-\varepsilon, \varepsilon[$ by $[0, \varepsilon[$ to obtain local semiflows if and only if the differential equation admits unique solutions around an initial value depending smoothly on the initial values. If for all $x \in U$ there is a local semiflow around $x$, then there is a maximal local semiflow. A maximal local semiflow $Fl : W \to E$ satisfies

$$(t, Fl(s, y)) \in W \text{ if and only if } (s + t, y) \in W$$

for all $(s, y) \in W$ and $t \in \mathbb{R}_{\geq 0}$. In this case $Fl(t, Fl(s, y)) = Fl(s + t, y)$. The notion of a semiflow is redundant on finite-dimensional vector spaces.

We are in particular interested in special types of differential equations on Fréchet spaces, namely Banach map perturbated bounded linear equations. Given a bounded linear operator $A : E \to E$, the abstract Cauchy problem associated to $A$ is given by the differential equation associated to $A$. We assume that there is a smooth semigroup of bounded linear operators $S : \mathbb{R}_{\geq 0} \times E \to E$ such that

$$\lim_{t \downarrow 0} \frac{S_t x - x}{t} = Ax$$

which is a global semiflow for the linear vector field $A$. Notice that the theory of bounded linear operators on Fréchet spaces contains as a special case Hille-Yosida-Theory of unbounded operators on Banach spaces (see for example [18]). Given a Banach map $P : U \subset E \to E$ with splitting $P = Q \circ R$ we want to investigate the solutions of

$$\frac{d}{dt} x(t) = Ax(t) + P(x(t))$$

We do this by constructing a solution to the integral equation arising from variation of constants

$$Fl(t, x) = S_t x + \int_0^t S_{t-s} P(Fl(s, x)) ds$$

for small positive time intervals and an open neighborhood of a given initial value. Given $x_0 \in U$ there exists - due to theorem 4 - a seminorm $p$ on $E$ and $\delta > 0$ such that

$$||R(x_1) - R(x_2)|| \leq p(x_1 - x_2)$$

for $p(x_i - x_0) < \delta$ and $i = 1, 2$, where $||.||$ denotes the norm on $X$. Furthermore given $y_0 \in X$, then for any seminorm $q$ on $F$ there are constants $C_q$ and $\delta_q$ such that

$$q(Q(y_1) - Q(y_2)) \leq C_q ||y_1 - y_2||$$

for $||y_i - y_0|| < \delta_q$ and $i = 1, 2$. By the uniform boundedness principle the set of continuous linear operators $\{S_t\}_{0 \leq t \leq T}$ is uniformly bounded for any fixed $T \geq 0$, i.e. for any seminorm $p$ on $E$ there is a seminorm $q_p$ such that

$$p(S_t x) \leq q_p(x)$$

for $t \leq T$. We denote by $C([0, \varepsilon], X)$ continuous curves on the interval $[0, \varepsilon]$ to $X$, $y_0 := R(x_0)$. Without any restriction we can assume that $x_0 = 0$ and $y_0 = 0$ by translations. We can then define a mapping

$$M : U' \subset E \times V' \subset C([0, \varepsilon], X) \to V'$$

such that $M(x, g)(t) = R(S_t x + \int_0^t S_{t-s} Q(g(s)) ds)$ for $t \in [0, \varepsilon]$. Given $g \in C([0, \varepsilon], X)$ such that $||g(t)|| \leq \theta$ for $0 \leq t \leq \varepsilon$ with $\{g | \sup_t ||g(t)|| \leq \theta\} \subset V$, we have

$$p(S_t x + \int_0^t S_{t-s} Q(g(s)) ds) \leq q_p(x) + \varepsilon(q_p(Q(y_0)) + C_{q_p}\theta)$$

provided $\theta \leq \delta_{q_p}$. This can be made smaller than $\theta$ if $\varepsilon$ is appropriately small and $U' := \{x \in E \text{ with } q_p(x) < \eta\}$ with $\eta$ appropriately small. In particular $C_{q_p}\varepsilon < 1$. If we assume these conditions, then $M$ is well defined, continuous and furthermore

$$\sup_t ||M(x, g_1)(t) - M(x, g_2)(t)|| \leq \varepsilon \sup_t q_p(Q(g_1(t)) - Q(g_2(t))) \leq$$
$$\leq C_{q_p}\varepsilon \sup_t ||g_1(t) - g_2(t)||$$

Consequently $M(x, .)$ is a contraction in $V'$ with contraction constant bounded uniformly in $x \in U'$ by a constant strictly smaller than 1. It follows that there is a unique $g(t, x)$ for any $x \in U'$ depending continuously on $x$, such that

$$M(x, g) = g$$

by the contraction mapping theorem. We define

$$Fl(t, x) := S_t x + \int_0^t S_{t-s} Q(g(s, x)) ds$$

and obtain

$$Fl(t, x) = S_t x + \int_0^t S_{t-s} P(Fl(s, x)) ds$$

since $R(Fl(t, x)) = g(t, x)$ by construction. By induction, smoothness with respect to time is easily established and uniqueness follows from the contraction mapping theorem. In the time dependent case - with respect to the bounded linear map and the Banach map - we can argue in the same manner and obtain a unique time-dependent semiflow, smooth with respect to time.

Concerning smoothness with respect to the initial value, we proceed in the following way: First we show that there exist directional derivatives: By Taylor's formula we obtain

$$P(x_1) - P(x_2) = L(x_1, x_0) \cdot (x_1 - x_2)$$

where $L(x_1, x_2) \cdot h := \int_0^1 DP(x_0 + s(x_1 - x_0)) \cdot h ds$ is a Banach map in all three variables. So we can solve the system given by

$$(x_0, x_1, h) \mapsto (Ax_0 + P(x_0), Ax_1 + P(x_1), Ah + L(x_1, x_2) \cdot h)$$

with a flow $Fl(t, x_0, x_1, h) = (Fl(t, x_0), Fl(t, x_1), M(t, x_0, x_1, h))$, smooth in time and continuous in initial values, where the dependence on $h$ is homogenous, so the flow can be defined everywhere in $h$. By uniqueness the identity

$$\frac{d}{dt}(Fl(t, x_0) - Fl(t, x_1))$$
$$= A(Fl(t, x_0) - Fl(t, x_1)) + L(Fl(t, x_0), Fl(t, x_1)) \cdot (Fl(t, x_0) - Fl(t, x_1))$$

leads to

$$M(t, x_0, x_1, x_0 - x_1) = Fl(t, x_0) - Fl(t, x_1).$$

By homogenity in $h$ we obtain the existence of the directional derivatives and its continuity in point and direction at the domain of definition. The argument follows an argument given in [9]. If the dependence were smooth, the derivative with respect to the initial value would satisfy the integral equations

$$D_2 Fl(t, x) \cdot v = S_t v + \int_0^t S_{t-s} \cdot DP(Fl(s, x)) \cdot D_2 Fl(t, x) \cdot v ds$$

$$Fl(t, x) = S_t x + \int_0^t S_{t-s} P(Fl(s, x)) ds$$

providing an equation of the same type with Banach mapping

$$(x, y) \mapsto (P(x), DP(x) \cdot y)$$

and semigroup $S_t \oplus S_t$ for $t \geq 0$. This equation has by the above procedure a unique solution, continuous with respect to initial values. If we integrate $u \mapsto D_2 Fl(t, x + uv) \cdot v$ with respect to $u$ we obtain $Fl(t, x + uv) - Fl(t, x)$ by uniqueness and Taylor's formula for $P$. Thus we can conclude by Theorem 12.8. of [12] and induction. There is one argument hidden, namely, that the domain in $t$ might shrink. We overcome this difficulty by the local existence result for the directional derivatives, which means that we can construct solutions for the derivative-equation by taking original solutions.

In an analogue manner we proceed in the time dependent case. Consequently we obtain smooth semiflows in all possible cases, which proves the following theorem:

**Theorem 7.** *Let $E$ be a Fréchet space and $I$ an interval, $A : I \times E \to E$ a smooth curve of bounded linear operators, such that there is a local time dependent semiflow $Ev(t, s, x)$ for $t \geq s \geq 0$ in appropriate domains and $x \in E$ satisfying*

$$\frac{d}{dt} Ev(t, s, x) = A(t) Ev(t, s, x)$$
$$Ev(t, t, x) = x$$

*for $t \geq s \geq 0$ and $x \in E$. If $P : I \times U \to E$ is a time-dependent family of Banach mappings, then there is a time dependent local semiflow $Fl(t, s, x)$ for $t \geq s \geq 0$ appropriate and $x \in E$ satisfying*

$$\frac{d}{dt} Fl(t, s, x) = A(t) Fl(t, s, x) + P(t, Fl(t, s, x))$$
$$Fl(t, t, x) = x$$

*for $t \geq s \geq 0$ and $x \in E$.*

**Remark 6.** *Parameter dependence can be easily treated by the above methods on the extended phase space.*

## 3. Submanifolds and Weak Foliations

We are interested in the geometry generated by a finite number of vector fields given on an open subset of a Fréchet space $E$. Therefore we need the notions of finite-dimensional submanifolds (with boundary) of a Fréchet space (see [12] for all details and more). First we introduce the notion of a smooth Fréchet manifold.

A chart on a set $M$ is a bijective mapping $u : U \to u(U) \subset E_U$, where $E_U$ is a Fréchet space and $U \subset M$, $u(U) \subset E_U$ is open. We shall denote a chart by $(U, u)$ or $(u, u(U))$. For two charts $(U_\alpha, u_\alpha)$, $(U_\beta, u_\beta)$ the chart changing are given by $u_{\alpha\beta} := u_\alpha \circ u_\beta^{-1} : u_\beta(U_{\alpha\beta}) \to u_\alpha(U_{\alpha\beta})$, where $U_{\alpha\beta} := U_\alpha \cap U_\beta$. An atlas is a collection of charts such that the $U_\alpha$ form a cover of $M$ and the chart changings are defined on open subsets of the respective Fréchet spaces. A $C^\infty$-atlas is an atlas with smooth chart changings. Two $C^\infty$-atlases are equivalent is their union is an $C^\infty$-atlas. A maximal $C^\infty$-atlas is called a $C^\infty$-structure on $M$ (maximal is understood with respect to some carefully chosen universe of sets). A (smooth) manifold is a set together with a $C^\infty$-structure. A topological manifold is a manifold with a $C^0$-structure. If we are given a Banach space, then we can define $C^k$-manifolds by an $C^k$-structure.

A smooth mapping $f : M \to N$ between smooth manifolds is defined in the canonical way, i.e. for any $x \in M$ there is a chart $(V, v)$ with $f(x) \in V$, a chart $(U, u)$ of $M$ with $x \in U$ and $f(U) \subset V$, such that $v \circ f \circ u^{-1}$ is smooth. This is the case if and only if $f \circ c$ is smooth for all smooth curves $c : \mathbb{R} \to M$, where the concept of a smooth curve is easily set upon.

The final topology with respect to smooth curves or equivalently the final topology with respect to all inverses of chart mappings is the canonical topology of the smooth manifold. We assume manifolds to be smoothly Hausdorff (see the discussion in [12], p. 265), i.e. the real valued smooth functions on $M$ separate points. The product of smooth manifolds is defined canonically by building up the product of the atlases.

A submanifold $N$ of a Fréchet manifold $M$ is given by a subset $N \subset M$, such that for each $n \in N$ there is a chart $(u, u(U))$, a splitting $E = E' \times E''$ and $u(U) = V \times W$ with $u(N) = V \times \{u(n)''\}$. By a splitting we shall always understand $E'$ and $E''$ as closed subspaces of $E$.

A finite-dimensional manifold with boundary of dimension $n$ is defined as ordinary manifold except that we take open subsets in a halfspace $\{x \in \mathbb{R}^n$ with $x_n \geq 0\}$. For the notion (without surprises) of smooth mappings on such open sets see any textbook on differential geometry, for example [13]. The boundary $\{x$ with $x_n = 0\}$ of the subspace models the boundary $\partial N$ of the manifold $N$, which is canonically a manifold without boundary of dimension $n - 1$. A submanifold with boundary is given by the analogue submanifold-structure.

We specialize to submanifolds with boundary of Fréchet spaces - the case we are interested in. Let $N \subset E$ be a submanifold with boundary, then the tangent space $T_n N$ at $n \in N$ is given by the derivatives $\frac{d}{dt}|_{t=0} c(t)$ of all smooth curves $c : \mathbb{R} \to N \subset E$ with $c(0) = n$. It is a vector space at any point $n \in N \setminus \partial N$ and a half space at $n \in \partial N$. Nevertheless we denote by $T_n N$ the generated vector space and by $(T_nN)_+$ the canonically given halfspace. There is a natural structure of a manifold with boundary on the disjoint union $\cup_{n\in N}T_nN$, the tangent bundle. A vector field on $N$ is a smooth mapping on $N$ associating to any point a tangent vector. A vector field $X$ on $U \subset E$ is therefore simply a smooth mapping $X : U \to E$, since the tangent space at $x \in U$ is $E$. The set of vector fields is denoted by $\mathfrak{X}(N)$. The derivative $Tf : TM \to TN$ of a smooth function $f : M \to N$ is defined in the following way

$$T_af(a, [c]) = (f(a), [f \circ c])$$

where $[c]$ denotes the equivalence class of curves $c : \mathbb{R} \to M \subset E$ having the same derivative at the $t = 0$. A local diffeomorphism $f : M \to N$ is an isomorphism on the respective tangent spaces. So the derivative mapping $T_xf : T_xU \to T_{f}(x)V$ is given by $T_xf(v) = Df(x)\dot{v}$.

Finite dimensional submanifolds (with boundary) of a Fréchet manifold $M$ can locally be given by immersions. An immersion is a smooth mapping $\phi : N \to M$ with $T_a\phi : T_aN \to T_{\phi(a)}M$ injective for all $a \in N$ and $N$ a finite-dimensional manifold of dimension $n$.

**Lemma 4** (Submanifolds by Parametrization). *Let $N$ be a finite-dimensional manifold (with boundary), $E$ a Fréchet space and $\phi : N \to E$ an immersion, then for any $n_0 \in N$ there is a small open neighborhood $V$ such that $\phi(V_1)$ is a submanifold of $E$.*

*Proof.* Given an immersion we shall construct a submanifold chart for the local image of $N$. We assume that $N$ is an open subset of $\mathbb{R}^m$ and - by translation - $\phi(n_0) = 0$, since it is a local result. Given a linear basis $e_1, ..., e_n$ of $T_{n_0}N$ we get linearly independent vectors $T_{n_0}(\phi)(e_i) =: f_i \in E$. We choose $l_1, ..., l_m$ linearly independent linear functionals, such that $l_i(f_j) = \delta_{ij}$ and get a splitting $E = E' \times E''$ with $\dim E' = m$ via $E'' := \cap_{i=1}^m \ker l_i$. The projection on the first variable $p_1$ induces a local diffeomorphism $p_1 \circ \phi$ on a small open neighborhood $V$ of $n_0 \in N$ by the classical inverse function theorem. The inverse is denoted by $\psi : V' \subset E' \to V$. Now we construct a new diffeomorphism

$$\eta(n, x'') = (p_1 \circ \phi(n), x'' + p_2 \circ \phi(n))$$

on $V \times W''$, which is invertible by the above considerations:

$$\eta^{-1}(y', y'') = (\psi(y'), y'' - p_2 \circ \phi(\psi(y'))),$$

$\eta^{-1}$ defines a submanifold chart for $\phi(V)$ since

$$\eta^{-1}(\phi(n)) = (n, 0)$$

for $n \in V$ by definition. The proof for manifolds with boundary works in the same way except that the linear basis at a boundary point $n \in \partial N$ has to lie in $(T_nN)_+$. □

Two vector fields $X \in \mathfrak{X}(M)$, $Y \in \mathfrak{X}(N)$ are called $f$-related for a smooth map $f : M \to N$ if $Tf \circ X = Y \circ f$. We obtain a bounded linear mapping $f^* : \mathfrak{X}(N) \to \mathfrak{X}(M)$, the pull back, for a local diffeomorphism $f : M \to N$ by the following formula

$$(f^*Y)(x) = (T_xf)^{-1}(Y_{f(x)})$$

for $Y \in \mathfrak{X}(N)$. The push forward for a diffeomorphism is defined by $f_* = (f^{-1})^*$.

Each vector field $X$ admitting a local flow has a maximal local flow $Fl^X$ producing all integral curves on a connected set, furthermore $X$ is $Fl_t^X$-related to itself for any $t$, since

$$T_x Fl_t^X(x) \cdot X(x) = \frac{d}{ds} Fl_t^X(Fl_s^X(x))|_{s=0} = \frac{d}{ds} Fl_{t+s}^X(x)|_{s=0} = X(Fl_t^X(x))$$

Given a local semiflow, then the solutions of the differential equations are unique: Given an integral curve $c$ we obtain

$$\frac{d}{ds} Fl_{t-s}(c(s)) = -X(Fl_{t-s}(c(s))) + X(Fl_{t-s}(c(s))) = 0$$
$$Fl_{t-s}(c(s)) = c(t)$$

for all $s \leq t$.

Given two vector fields $Y \in \mathfrak{X}(M)$ and $X \in \mathfrak{X}(M)$ admitting a local flow, where $M$ either denotes a finite-dimensional submanifold of a Fréchet space or an open subset, we can define the Lie bracket, the most important notion of differential geometry:

$$[X, Y] = \frac{d}{dt}(Fl_{-t}^X)^* Y|_{t=0}$$

Notice furthermore that for any local diffeomorphism $\phi$

$$\phi^*[X, Y] = [\phi^* X, \phi^* Y]$$

so the pull back is a bounded Lie algebra homomorphism, since vector fields constitute a Lie algebra with the Lie bracket. For details in all possible directions see [12].

Finally we provide a formula for the Lie bracket on a Fréchet space. Given two vector fields $X$ and $Y$ on an open subset $U \subset E$, where $X$ admits a flow $Fl^X$, then

$$\begin{aligned}
[X, Y](x) &= \frac{d}{dt}(Fl_{-t}^X)^* Y(x)|_{t=0} \\
&= \frac{d}{dt} DFl_t^X(Fl_{-t}^X(x)) \cdot Y(Fl_{-t}^X(x))|_{t=0} \\
&= DX(x) \cdot Y(x) - D^2 Fl_0^X(x)(X(x), Y(x)) - DY(x) \cdot X(x) \\
&= DX(x) \cdot Y(x) - DY(x) \cdot X(x)
\end{aligned}$$

for $x \in U \subset E$ (see [12] and [11]). Due to this formula we can immediately conclude some properties of the submodule of Banach map vector fields:

**Lemma 5.** *Let $U$ be an open set in a Fréchet space $E$, then $\mathcal{B}(U)$ is a subalgebra with respect to the Lie bracket. Let $A$ be a bounded linear operator on $E$, then $[A, \mathcal{B}(U)] \subset \mathcal{B}(U)$. Consequently the Lie algebra $L(E)$ acts on $\mathcal{B}(U)$ by the Lie bracket.*

*Proof.* Given two Banach maps $P_1$ and $P_2$, $DP_1(x) \cdot P_2(x) = DQ_1(R_1(x)) \cdot DR_1(x) \cdot P_2(x)$ holds, which can be written as composition of $DQ_1(v) \cdot w$ for $v, w \in X$ and $(R_1(x), DR_1(x) \cdot P_2(x))$ for $x \in U$. So the Lie bracket lies in $\mathcal{B}(U)$. Given $A \in L(E)$, we see that $AP_1(x) - DP_1(x) \cdot Ax$ is a Banach map by an obvious decomposition. ☐

**Definition 4.** *Let $E$ be a Fréchet space, $U$ an open subset. A distribution on $U$ is a collection of vector subspaces $D = \{D_x\}_{x \in U}$ of $\{T_x U\}_{x \in U}$. A distribution is* *called smooth if for any $x \in U$ there is an open neighborhood $x \in V \subset U$ and $n$ vector fields $X_1, ..., X_n$ with*

$$\langle X_1(y), ..., X_n(y) \rangle = D_y$$

*for $y \in V$. A distribution $D$ on $U$ is said to be involutive if for any two locally given vector fields $X, Y$ with values in $D$ the Lie bracket $[X, Y]$ has values in $D$. A distribution is said to have constant rank if $\dim_{\mathbb{R}} D_x$ is locally constant on $U$.*

We denote by $\langle \ldots \rangle$ the generated vector space over the reals $\mathbb{R}$. In the case of a smooth distribution $D$ we sometimes apply the notation $D = \langle X_1, \ldots, X_n \rangle$ which means that $D_x$ is vector space generated by $X_1(x), \ldots, X_n(x)$. We call this the linear span of $X_1, \ldots, X_n$. If the vector fields $X_1, \ldots, X_n$ are linearly independent in a open neighborhood of a point and if they span $D$, then we call these vector fields a local basis.

**Remark 7.** *Given a smooth distribution $D$ on $U$ such that the dimensions of $D_x$ are bounded by a fixed constant $N$. Let $x \in U$ be a point with maximal dimension $n_x = \dim_{\mathbb{R}} D_x$, then there are $n$ smooth vector fields $X_1, ..., X_n$ and an open subset $x \in V \subset U$ such that*

$$\langle X_1(y), ..., X_n(y) \rangle = D_y$$

*for $y \in V$. At $x$ there are $n_x$ linearly independent vectors, say $X_1(x), ..., X_{n_x}(x)$. Choosing $n_x$ continuous linear functionals $l_1, ..., l_{n_x} \in E'$ with $l_i(X_j(x)) = \delta_{ij}$, then the continuous mapping $\mathbf{l} : U \to L(\mathbb{R}^{n_x})$, $x \mapsto (l_i(X_j(x)))$ has range in the invertible matrices in a small neighborhood of $x$. Consequently in this neighborhood the dimension of $D_y$ is at least $n_x$. It follows by maximality of $n_x$ that it is constant. In particular a smooth distribution with constant rank admits locally a basis.*

The concept of weak foliations will be perfectly adapted to our purposes in HJM-theory:

**Definition 5.** *Given a distribution $D$ on an open subset $U \subset E$, a tangential manifold (with boundary) for $D$ is an immersion $\phi_N : N \to U$ such that $D_n \subset T_n \phi_N(T_n N)$ for all $n \in N$. We say that $D$ admits a weak foliation of rank $m$ if for any point $x_0 \in U$ there is an open subset $V$ of $\{\mathbf{u} \in \mathbb{R}^m \text{ with } u_m \geq 0\}$, an open subset $W''$ of a Fréchet space $E''$ and a smooth mapping $\psi : V \times W'' \to U$ such that*

i) *There is a point in the parameter space such that $x_0 = \psi(\mathbf{u}_0, x_0'')$;*
ii) *$T_{(\mathbf{u}_0, x'')}\psi : \mathbb{R}^m \times E'' \to E$ is an isomorphism;*
iii) *$\psi(., x'')$ is an tangential manifold (with boundary) for $D$ for any $x'' \in W''$;*
iv) *The tangent distribution given by all tangent spaces of the tangential manifolds has constant rank $m$.*

**Remark 8.** *The immersed manifolds are called leafs. In the chart sometimes the notion plaque is applied for their local image. The dimension $m$ of the tangential manifolds will generically be greater than the dimension of $D_x$ of the given distribution $D$.*

Classically one is interested in the existence of weak foliations for a given distribution of minimal dimension $m$. Therefore we shall need the following essential lemma:

**Lemma 6.** *Let $D$ be an involutive, smooth distribution of constant rank $n$ on an open subset $U$ of a Fréchet space $E$. Let $X$ and $Y$ be vector fields with values in $D$ and let $X$ admit a local flow, then*

$$(Fl_t^X)^*(Y)(x) \in D_x$$

*for $x \in U$, where it is defined.*

*Proof.* Given a local basis $X_1, ..., X_n$ around $x_0$, we have by involutivity that $[X, X_i] = \sum_{k=1}^{n} f_i^k X_k$. Remark that $f_i^k$ are smooth functions locally around $x_0$: Given $n$ linear independent functionals $l_m$ such that $l_m(X_j(x_0)) = \delta_{mj}$, then

$$l_m([X, X_i](x)) = \sum_{k=1}^{n} f_i^k l_m(X_k(x))$$

Since the matrix $M(x) := (l_m(X_k(x)))$ is invertible at $x_0$ and has smooth entries, it is invertible on an open neighborhood of $x_0$, and the inverse has smooth entries. The smooth inverse matrix applied to the left hand vector proves the smoothness of $f_i^k$. By involutivity we obtain

$$\frac{d}{dt}(Fl_t^X)^*(X_i) = -\frac{d}{ds}(Fl_{t-s}^X)^*(X_i)|_{s=0} =$$

$$= -\frac{d}{ds}(Fl_{-s}^X)^*(Fl_t^X)^*(X_i)|_{s=0} =$$

$$= -[X, (Fl_t^X)^*(X_i)] =$$

$$= -(Fl_t^X)^*[X, X_i] =$$

$$= -\sum_{k=1}^{n} f_i^k \circ Fl_t^X (Fl_t^X)^*(X_k)$$

which is a linear equation with time-dependent coefficients $g_i^k(t) := -f_i^k(Fl_t^X(x))$ for $(Fl_t^X)^*(X_i)(x)$ at any point $x$ in that neighborhood of $x_0$. Consequently provided the initial values lie in $D_x$ the solution lies in $D_x$, but $(Fl_0^X)^*(X_i)(x) = X_i(x)$.  □

Now we can state and prove the Frobenius theorem for weak foliations, remark that in this infinite-dimensional case the Frobenius condition is a sufficient condition for the existence of a weak foliation, since there are problems to solve ordinary differential equations. Nevertheless we shall see that there is a partial answer in the other direction.

**Theorem 8.** *Let $D$ be an involutive, smooth distribution of constant rank $n$ on an open subset $U$ of a Fréchet space $E$ and assume that for any point $x_0$ there is a basis of the distribution $X_1, ..., X_n$ of locally defined vector fields, where $X_1, ..., X_{n-1}$ admit local flows $Fl_t^{X_i}$ and $X_n$ admits a local semiflow, then $D$ admits a weak foliation of rank $n$.*

*Proof.* There are locally on some open set $n$ linearly independent vector fields $X_1, ..., X_n$ generating each $D_x$ with local flows $Fl^{X_i}$ for $i = 1, ..., n-1$ and a local semiflow $Fl^{X_n}$. Without restriction we set $x_0 = 0$. We choose $n$ linear functionals $l_1, ..., l_n$ such that $l_i(X_j(0)) = \delta_{ij}$, so we get some splitting $r : E \to \mathbb{R}^n \times E''$. We define the parametrization $\alpha(\mathbf{u}, y) = Fl_{u_1}^{X_1} \circ ... \circ Fl_{u_n}^{X_n}(0, y)$ for $(\mathbf{u}, y)$ in an open subset around $0$ in $\mathbb{R}^n \times E''$. We can calculate by looking at the partial derivatives the tangent map applied to the canonical standard basis $\frac{\partial}{\partial u_i}|_{\mathbf{u}} := [s \mapsto \mathbf{u} + se_i]$ of $(T_{\mathbf{u}}\mathbb{R}^n)_+$, the

$$T_{\mathbf{u}}\alpha(\frac{\partial}{\partial u_1}|_{\mathbf{u}}, ..., \frac{\partial}{\partial u_n}|_{\mathbf{u}})(X_1, ..., (Fl_{u_1}^{X_1})_* ... (Fl_{u_{n-1}}^{X_{n-1}})_* X_n)(\alpha(\mathbf{u}, y))$$

which lie by Lemma 6 in $D_{\alpha(\mathbf{u},y)}$ and generate it in a small neighborhood by a dimension argument. It is essential that the first $n-1$ vector fields admit a local flow. So we obtain a family of tangential manifolds for $D$. Each parametrization for fixed $y$ defines locally a smooth submanifold with boundary $\alpha(u_1, .., u_{n-1}, 0, y)$. □

**Remark 9.** *For details on Frobenius theorems in the classical setting see [11]. The phenomenon that there is no Frobenius chart is due to the fact that there is one vector field among the vector fields $X_1,...,X_n$ (generating the distribution $D$) admitting only a local semiflow. If all of them admitted flows, there would exist a Frobenius chart, which can be given by a construction outlined in [19]. The non-existence of a Frobenius-chart means that the leafs cannot be parallelized, since they follow semiflows, which means that "gaps" between two leafs can appear. These "gaps" do not appear as interior point of a leaf. This is an infinite dimensional phenomenon, which does not appear in finite dimensions.*

**Remark 10.** *Let $D$ be a smooth distribution of constant rank $n$ on an open subset $U$ of a Fréchet space $E$ and assume that for any point $x_0$ there is a basis of the distribution $X_1,...,X_n$ of locally defined vector fields, where $X_1,...,X_{n-1}$ admit local flows $Fl_t^{X_i}$ and $X_n$ admits a local semiflow. Assume furthermore that $D$ admits a weak foliation of dimension $n$, then the distribution is involutive. This is easily seen by applying $(Fl_t^{X_i})^*$ to $X_j$ at $x \in U$, which is by assumption an element of $D_x$ (since the flow restricts by uniqueness of integral curves to the tangential manifolds), so differentiation and smoothness yield the result by the formula of lemma 6.*

## 4. Existence of Finite dimensional Realizations

Subsequently, we let assumptions (A1)–(A4) be in force. For simplicity we shall discuss only the time-homogeneous case (see Remark 14 for an outline of the time-inhomogeneous case).

The vector fields $\{\mu, \sigma^1, \ldots, \sigma^d\}$ induce two smooth distributions on $D(A^\infty)$: their linear span $D$, and $D_{LA}$, the linear span of all multiple Lie brackets of these vector fields.

Recapturing the discussion in Section 1 we now can say that (1.1) has generically an FDR if $D$ admits locally a weak foliation. We shall see below (Proposition 2) that the assumption that $D$ admits a weak foliation on some open set $V$ is equivalent to the assumption that $D_{LA}$ is an involutive, smooth distribution with constant rank on some open (in general smaller) subset $U \subset V$.

Define the smooth map $\Gamma := \sum_{j=1}^{d} \Gamma^j : B \to D(A^\infty)$ by

$$\Gamma^j(y) := \mathcal{S}(\phi^j(y)) - \frac{1}{2}D\phi^j(y)\left(\ell(\phi^j(y))\right).$$

So that we can write $\mu(h) = Ah + \Gamma(\ell(h))$.

Let $U$ denote an open set in $D(A^\infty)$ in what follows.

**Lemma 7.** *Let $X_1, \ldots, X_k$ be linearly independent Banach-maps on $U$, for some $k \in \mathbb{N}$. Then the set*

$$\mathcal{N} = \{h \in U \mid \mu(h) \in \langle X_1(h), \ldots, X_d(h) \rangle\}$$

*is nowhere dense in $U$.*

*Proof.* We argue by contradiction. Suppose there exists a set $V \subset \overline{\mathcal{N}}$ which is open in $D(A^\infty)$. For any $h \in V \cap \mathcal{N}$, there exist numbers $c_1(h), \ldots, c_d(h)$ such that

$$\mu(h) = \sum_{j=1}^{d} c_j(h) \sigma^j(h). \tag{4.1}$$

Since $D(A^*)$ is dense in $H$, we can find $\xi_1, \ldots, \xi_d \in D(A^*)$ such that the $d \times d$-matrix $M_{ij}(h) := \xi_i(\sigma^j(h))$ is smooth and invertible on $V$ (otherwise we choose a smaller open subset $V$). Hence

$$\begin{pmatrix} c_1(h) \\ \vdots \\ c_d(h) \end{pmatrix} = M^{-1}(h) \begin{pmatrix} \xi_1(\mu(h)) \\ \vdots \\ \xi_d(\mu(h)) \end{pmatrix}$$

can be smoothly extended to $V$. Then (4.1) implies that $A$ is a Banach-map on $V$. But this contradicts (A4), whence the claim. $\qquad \square$

In what follows, we suppose that $D_{LA}$ has constant dimension $k_D \in \mathbb{N}$ on $U$.

**Lemma 8.** *We have $k_D > d$ and*

$$\mu(h) \notin \langle \sigma^1(h), \ldots, \sigma^d(h) \rangle, \quad \forall h \in U. \tag{4.2}$$

*Moreover, for any $h_0 \in U$ there exists an open neighborhood $V$ and Banach-maps $X^{d+1}, \ldots, X^{k_D-1}$ on $V$ such that $\{\mu, \sigma^1, \ldots, \sigma^d, X^{d+1}, \ldots, X^{k_D-1}\}$ is a basis for $D_{LA}$ on $V$.*

*Proof.* Assumption (A3) implies $k_D \geq d$. By Lemma 7 there exists a nowhere dense set $\mathcal{N} \subset U$ such that $\mu(h) \notin \langle \sigma^1(h), \ldots, \sigma^d(h) \rangle$, for all $h \in U \setminus \mathcal{N}$. Therefore $k_D > d$.

Now suppose $\mu(h_0) \in \langle \sigma^1(h_0), \ldots, \sigma^d(h_0) \rangle$, for some $h_0 \in U$. By the definition of $D_{LA}$ and Lemma 5 there exist $k_D - d$ Banach-maps $X^{d+1}, \ldots, X^{k_D}$ on $U$ such that

$$D_{LA}(h) = \langle \sigma^1(h), \ldots, \sigma^d(h), X^{d+1}(h), \ldots, X^{k_D}(h) \rangle,$$

for $h = h_0$, and hence for all $h$ in a neighborhood of $h_0$, by continuity. But this implies that $\mu(h)$ lies in the span of Banach-maps, for all $h \in V$. This contradicts Lemma 7, whence the claim follows. $\qquad \square$

**Remark 11.** *Lemma 8 proves a conjecture in [3], namely that any short rate realization is of dimension 2 (see Remark 7.1 therein).*

In view of Lemma 8 the minimal generic FDRs are of dimension $d+1$. The following result is an extension of [3, Proposition 7.4].

**Proposition 1.** *Suppose $k_D = d + 1$. Then there exist $d$ linear independent (constant) vectors $\Lambda_1, \ldots, \Lambda_d \in D(A^\infty)$ such that*

$$\sigma^j(h) \in \langle \Lambda_1, \ldots, \Lambda_d \rangle, \quad \forall h \in U, \quad j = 1, \ldots, d.$$

Hence *necessarily* we are in the (slightly modified) "deterministic direction volatility" case of [3] (see Section 6 therein).

*Proof.* By (the proof of) Lemma 8 we have, for any $k = 1, \ldots, d$,

$$[\mu, \sigma^k](h) = \sum_{j=1}^{d} c_j^k(h)\sigma^j(h), \quad \forall h \in U.$$

Fix $h_0 \in U$. As in the proof of Lemma 7, we find an open neighborhood $V$ of $h_0$ and $\xi_1, \ldots, \xi_d \in D(A^*)$ such that the $d \times d$-matrix $M_{ij}(y) := \xi_i(\phi^j(y))$ is smooth and invertible for all $y \in W := \ell(V)$. In view of (A2), the sets $W$ and $W' := (\ell, \ell \circ A)(V)$ are open in $B$ and $B \times B$, respectively. A calculation shows that

$$[\mu, \sigma^k](h) = \Delta^k(\ell(h), \ell(Ah)),$$

where

$$\Delta^k(y, z) := A\phi^k(y) + D\Gamma(y)\ell(\phi^k(y)) - D\phi^k(y)\left(z + \ell(\Gamma(y))\right)$$

(linearity of $\ell$ is essential for obtaining this implicit dependence on $h$). Consequently, the functions

$$\gamma_j^k(y, z) := \sum_{i=1}^{d} M_{ji}^{-1}(y)\xi_i\left(\Delta^k(y, z)\right) : W' \to \mathbb{R}$$

are smooth and satisfy

$$\Delta^k(y, z) = \sum_{j=1}^{d} \gamma_j^k(y, z)\phi^j(y), \quad \forall (y, z) \in W'. \tag{4.3}$$

Differentiation of (4.3) with respect to $z$ (which makes sense since $W'$ is open) yields

$$D\phi^k(y) = \sum_{j=1}^{d} D_z \gamma_j^k(y, z)\phi^j(y), \quad \forall (y, z) \in W'.$$

Arguing again by linear independence of $\{\phi^1, \ldots, \phi^d\}$, we see that

$$D_z \gamma_j^k(y, z) \equiv: \beta_j^k(y)$$

are smooth functions of $y$ only, with values in the dual space of $B$. Hence $\phi$ satisfies a linear differential equation on $W$,

$$D\phi^k(y) = \sum_{j=1}^{d} \beta_j^k(y)\phi^j(y).$$

We may assume that $W$ is star-shaped with respect to $y_0 = \ell(h_0)$; that is,

$$y_0 + t(y - y_0) \in W, \quad \forall t \in [0, 1], \quad \forall y \in W.$$

Otherwise we replace $V$ by $\ell^{-1}(B_\varepsilon(y_0))$, where $B_\varepsilon(y_0) := \{y \in B \mid \|y - y_0\|_B < \varepsilon\}$, for some $\varepsilon > 0$ small enough. Let $y \in W$ and define $\psi^k(t) := \phi^k(y_0 + t(y - y_0))$. Then there exists an open interval $I$, which contains $[0, 1]$ and such that

$$\begin{cases} \partial_t \psi^k(t) = \sum_{j=1}^{d} \beta_j^k(y_0 + t(y - y_0))(y - y_0)\psi^j(t) \\ \psi^k(0) = \phi^k(y_0), \quad k = 1, \ldots, d, \end{cases}$$

for $t \in I$. This differential equation has a unique solution, which is of the form

$$\psi^k(t) = \sum_{j=1}^{d} \alpha_j^k(t) \phi^j(y_0).$$

In particular, $\phi^k(y) = \psi^k(1) \sum_{j=1}^{d} \alpha_j^k(1) \phi^j(y_0)$.

We have thus shown that, for any $y \in W$, there exists a $d \times d$-matrix $A_j^k(y)$ such that

$$\phi^k(y) = \sum_{j=1}^{d} A_j^k(y) \phi^j(y_0). \tag{4.4}$$

By linear independence of the $\phi^k$s we conclude that $(A_j^k)$ is smooth and invertible on $W$. Now the proposition follows by a continuity argument. ∎

Let the assumptions of Proposition 1 be in force. We now shall see how an FDR looks. Note that we have

$$D\sigma^j(h)\sigma^j(h) \in \langle \Lambda_1, \dots, \Lambda_d \rangle, \quad \forall h \in U. \tag{4.5}$$

Write

$$\nu(h) := Ah + \alpha_{HJM}(h).$$

By (4.5) it follows that $D_{LA} = \langle \nu, \Lambda_1, \dots, \Lambda_d \rangle_{LA} = \langle \nu, \Lambda_1, \dots, \Lambda_d \rangle$, where the latter equality can be deduced as in the proof of Lemma 8. Also we see that

$$[\nu, \Lambda_k](h) = D\nu(h)\Lambda_k \in \langle \Lambda_1, \dots, \Lambda_d \rangle, \tag{4.6}$$

for all $h \in U$. Now let $h_0 \in U$ and $\mathcal{M}$ be a leaf of the weak foliation of $D$ in $U$ through $h_0$ (and hence a $(d+1)$-dimensional tangential manifold for $\langle \nu, \Lambda_1, \dots, \Lambda_d \rangle$). As in the proof of Theorem 8 we obtain a parametrization of $\mathcal{M}$ by

$$\alpha(u, y) = Fl_u(h_0) + \sum_{k=1}^{d} y_k \Lambda_k, \quad (u, y) \in [0, \varepsilon) \times V, \tag{4.7}$$

for some open set $V \subset \mathbb{R}^d$ and $\varepsilon > 0$, where $Fl_u$ is the smooth semiflow induced by $\nu$.

Though the FDRs are already nicely specified by the parametrization (4.7), the following geometric picture is worth noticing. Let $\mathcal{M}$ be a leaf through $h_0 \in U$ as above. The band of $d$-dimensional affine submanifolds of $\mathcal{M}$,

$$\mathcal{N}(u) := (Fl_u(h_0) + \langle \Lambda_1, \dots, \Lambda_d \rangle) \cap \mathcal{M}, \quad u \in [0, \varepsilon),$$

forms a foliation of $\mathcal{M}$. Using Taylor's formula we calculate, for $h \in \mathcal{M}$,

$$\nu(h) = \nu(Fl_u(h_0)) + \int_0^1 (D\nu(Fl_u(h_0) + s(h - Fl_u(h_0)))(h - Fl_u(h_0))) \, ds$$
$$=: \nu(Fl_u(h_0)) + \tilde{\nu}(u, h).$$

By (4.6) we have $\tilde{\nu}(u, h) \in \langle \Lambda_1, \dots, \Lambda_d \rangle$, for all $h \in \mathcal{N}(u)$, and we already know that $\nu(Fl_u(h_0)) \notin \langle \Lambda_1, \dots, \Lambda_d \rangle$. In other words, the vector-field $\nu$ restricted to $\mathcal{N}(u)$ splits into a component lying in $\langle \Lambda_1, \dots, \Lambda_d \rangle$ and a constant "outward pointing" component.

It remains to see whether more can be said about the functions $\Lambda_1, \ldots, \Lambda_d$. We can write $\sigma^i(h) = \sum_{j=1}^d \rho^{ij}(h)\Lambda_j$, for some smooth invertible matrix-function $\rho(h) = (\rho^{ij}(h))$. Accordingly

$$\nu(h) = Ah + \sum_{i,j=1}^d a^{ij}(h)\mathcal{S}(\Lambda_i, \Lambda_j),$$

where $a^{ij}(h) := (\rho\rho^*)^{ij}(h) = \sum_{l=1}^k \rho^{li}(h)\rho^{lj}(h)$. Now (4.6) implies

$$A\Lambda_k + \sum_{i,j=1}^d \left(Da^{ij}(h)\Lambda_k\right)\mathcal{S}\left(\Lambda_i, \Lambda_j\right) = \sum_{i=1}^d \gamma^{ki}(h)\Lambda_i.$$

Expressed as a point-wise equality of functions this reads

$$\partial_x \Lambda_k(x) + \sum_{i,j=1}^d \Gamma^{k,ij}(h)\partial_x\left(\Delta_i(x)\Delta_j(x)\right) = \sum_{i=1}^d \gamma^{ki}(h)\Lambda_i(x), \qquad (4.8)$$

where $\Delta_i(x) := \int_0^x \Lambda_i(\eta)\, d\eta$ and $\Gamma^{k,ij}(h) := (1/2)Da^{ij}(h)\Lambda_k$. Integrating with respect to $x$ yields

$$\partial_x \Delta_k(x) = \Lambda_k(0) + \sum_{i=1}^d \gamma^{ki}(h)\Delta_i(x) - \sum_{i,j=1}^d \Gamma^{k,ij}(h)\Delta_i(x)\Delta_j(x), \qquad (4.9)$$

and this has to hold for any $h \in U$. Thus any $h \in U$ implies a system of Riccati equations for the functions $\Delta_1, \ldots, \Delta_d$ (which have to hold simultaneously).

Suppose for the moment that the functions

$$\Delta_1, \ldots, \Delta_d \text{ and } \Delta_i\Delta_j, \ 1 \le i \le j \le d, \text{ are linearly independent.} \qquad (4.10)$$

Then (4.9) implies that $\gamma^{ki}(h) \equiv \gamma^{ki}$ and $\Gamma^{k,ij}(h) \equiv \Gamma^{k,ij}$ are constant on $U$. Choose continuous linear functionals $v_{ij}$ on $D(A^\infty)$ (they can be chosen from $H$) such that

$$v^{ij}(\Lambda_k) = \Gamma^{k,ij}.$$

Write $\delta^{ij}(h) := a^{ij}(h) - v^{ij}(h)$. Then we have

$$D\delta^{ij}(h)\Lambda_k = 0.$$

Hence $\delta^{ij}(h) = \delta^{ij}([h])$ is a function of the equivalence class $[h]$ of $h \in U$ given by the equivalence relation $f \sim g :\Leftrightarrow f - g \in \langle \Lambda_1, \ldots, \Lambda_d \rangle$. Thus

$$a^{ij}(h) = \delta^{ij}([h]) + v^{ij}(h). \qquad (4.11)$$

The parametrization (4.7) can always be globally extended for $(u, y) \in \mathbb{R}_{\ge 0} \times \mathbb{R}^d$ (though the image manifold can accumulate on itself). This does not mean, however, that a generic $(d+1)$-dimensional realization for (1.1) exists globally. As an illustration we remark that $a = (a^{ij})$ is by construction a strictly positive definite symmetric $d \times d$-matrix. Thus (4.11) yields the constraints

$$v^{ii}(h) + \delta^{ii}([h]) > 0, \quad i = 1, \ldots, d,$$

which in general are not globally satisfied. For a concrete example see the following section.

Let us finally consider the particular case where $\rho(h) \equiv \rho$ (and hence $\sigma^1, \ldots, \sigma^d$ are constant). Here (4.8) or (4.9), respectively, simplify to

$$\partial_x \Lambda_k(x) = \sum_{i=1}^d \gamma^{ki}(h) \Lambda_i(x). \tag{4.12}$$

Since $\Lambda_1, \ldots, \Lambda_d$ are linearly independent this implies that $\gamma^{ki}(h) \equiv \gamma^{ki}$ are constant. Equation (4.12) is a system of linear ODEs for the functions $\Lambda_i$.

A necessary condition for constant dimension, $k_D = d + 1$, of $D_{LA}$ on $U$ is (4.2). This yields

$$\nu(h) \notin \langle \Lambda_1, \ldots, \Lambda_d \rangle, \quad \forall h \in U. \tag{4.13}$$

Integrating with respect to $x$ it is easy to see that $\nu(h) \in \langle \Lambda_1, \ldots, \Lambda_d \rangle$ if and only if

$$h - h(0) + \frac{1}{2} \sum_{i,j=1}^d a^{ij} \Delta_i \Delta_j \in \langle \Delta_1, \ldots, \Delta_d \rangle. \tag{4.14}$$

This again holds if and only if $h \in \mathcal{O}$ where

$$\mathcal{O} := \left\{ b\mathbf{1} + \sum_{j=1}^d c_j \Delta_j - \frac{1}{2} \sum_{i,j=1}^d a^{ij} \Delta_i \Delta_j \mid b, c_1, \ldots, c_d \in \mathbb{R} \right\}$$

($\mathbf{1}$ denotes the constant function $\mathbf{1}(x) \equiv 1$). We have $\dim D_{LA} = d$ on $\mathcal{O}$. Nevertheless, $\mathcal{O}$ is a $(d+1)$-dimensional (affine) submanifold[2] of $D(A^\infty)$. Hence, although $D_{LA}$ does not have full rank, $d + 1$, globally on $D(A^\infty)$ we obtain a global generic $(d+1)$-dimensional realization. That is, a global weak foliation of $D(A^\infty)$ with leafs given by (4.7) and $\mathcal{O}$, respectively. The global HJM-model is by construction Gaussian.

**Remark 12.** *The ODEs (4.9) and (4.12), respectively, have already been derived and discussed in [3]. They are further analyzed in [2], where also the dynamics of the coordinate process is derived.*

**Remark 13.** *We can interpret the coordinates of the invariant submanifolds economically. If we turn back to the construction of the charts as given in Theorem 8 we see that we have been choosing continuous linear functionals $l_1, \ldots, l_n$. They can be chosen as linear combinations of point evaluations (we assumed them to be continuous) $ev_{x_1}, \ldots, ev_{x_n}$ by Theorem 23.10 in [12], which tells that linear combinations of point evaluations are dense in $C^\infty(\mathbb{R}_{\geq 0}, \mathbb{R})'$. For the proof look at the bounded multilinear mapping $(ev_{x_1}, \ldots, ev_{x_n}) \mapsto \det(ev_{x_i}(X_j(y)))$ which cannot be identically zero, since then the extension of the multilinear mapping to $(l_1, \ldots, l_n) \mapsto \det(l_i(X_j(y)))$ would be identically zero, too. However, this is impossible due to linear independence of $X_j(y)$. Consequently we can interpret the finite dimensional realization as Markovian process in "benchmarks". Notice that even though "benchmark" are obtained by linear projections on the forward rate, the Markovian process of the finite dimensional realization is in general not the projection of the solution of the HJM-equation.*

---

[2] We assume here that $\Delta_i, \Delta_i \Delta_j \in D(A^\infty)$. This requires a nice choice of $\gamma^{ki}$ in (4.12). That is, the matrix $(\gamma^{ij})$ must have negative eigenvalues.

The necessity of the assumptions in Proposition 1 for the existence of generic FDRs is essentially given by

**Proposition 2.** *Suppose that $D$ is a smooth distribution of constant rank $d + 1$ on an open subset $U$. Assume furthermore that on $U$ the distribution $D$ admits a weak foliation of dimension $m$ (by Definition 5 necessarily $m \geq d + 1$). Then $\dim_{\mathbb{R}} (D_{LA})_x \leq m$ and there are points $x_0 \in U$ such that $D_{LA}$ is a smooth involutive distribution of constant rank in a neighborhood of $x_0$.*

*Proof.* We repeat the argument of Remark 10. Take two vector fields of the canonical basis $\{\mu, \sigma^1, \dots, \sigma^d\}$ of $D$ . They admit a local (semi)flow, which restricts necessarily to the leafs of the weak foliation by assumption, so the formula of Lemma 6 tells that their Lie bracket lies in the tangent distribution of the weak foliation. By the special choice of the vector fields we know from section 2 that the Lie brackets admit local flows as Banach maps, so by uniqueness of the integral curve the flow restricts to the leafs. We can prove now by induction that $D_{LA}$ is an subdistribution of the tangent distribution of the weak foliation. Since the dimensions are globally bounded by $m$, it is a smooth distribution and therefore involutive around some points where the rank is maximal. $\qquad \square$

The general case where $D_{LA}$ has locally constant rank $k_D$ is sketched in the following proposition, which provides the full classification of FDRs in our framework. The detailed analysis and the geometric implications will be discussed elsewhere.

**Proposition 3.** *Suppose that $D_{LA}$ has locally constant dimension $k_D$ and that the mappings*

$$(l, l \circ A, \dots, l \circ A^r) : D(A^\infty) \to B^{r+1}$$

*are open for all integers $r \geq 0$. Then there exist $k_D - 1$ linearly independent vectors $\Lambda_1, \dots, \Lambda_{k_D-1} \in D(A^\infty)$ such that*

$$D_{LA} = \langle \mu, \Lambda_1, \dots, \Lambda_{k_D-1} \rangle$$
$$\sigma^j(h) \in \langle \Lambda_1, \dots, \Lambda_{k_D-1} \rangle$$

*locally in $h$.*

*Proof.* We proceed as in the proof of Proposition 1. Fixing a point $h_0$, we can choose linearly independent Banach map vector fields $X_1, \dots, X_{k_D-1}$ and an integer $r \geq -1$ such that

$$X^i(h) = \phi^i \circ (l, \dots, l \circ A^r)(h)$$

with $\phi^i : V \times B^r \subset B^{r+1} \to D(A^\infty)$ smooth with $V$ simply connected and

$$D_{LA} = \langle \mu, X_1, \dots, X_{k_D-1} \rangle$$

locally (a representation by constant vector fields is encoded with $r = -1$). We demand $r$ to be the *minimal* integer with the above properties. That is, either $r \geq 0$ and $D_{z^r} \phi^i \not\equiv 0$, for some $i$, or $r = -1$. We have to show that the latter is true.

We argue by contradiction and suppose that $r \geq 0$. Notice that $X_1, \dots, X_{k_D-1}$ spans exactly the subdistribution of $D_{LA}$ generated by Banach map vector fields. Consequently we obtain $[\mu, X^k](h) = \sum_{j=1}^{k_D-1} c_j^k(h) X^j(h)$ locally in $h$, since the Lie bracket is a Banach map. We can write

$$[\mu, X^k](h) = \Delta^k(l(h), \dots, l \circ A^{r+1}(h))$$

with

$$\Delta^k(z^0,...,z^{r+1}) = A\phi^k(z^0,...,z^r) + D\Gamma(z^0) \cdot l(\phi^k(z^0,...,z^r)) -$$
$$- D\phi^k(z^0,...,z^r) \cdot \left[ \begin{pmatrix} z^1 \\ \vdots \\ z^{r+1} \end{pmatrix} + \begin{pmatrix} l(\Gamma(z^0)) \\ \vdots \\ (l \circ A^r)(\Gamma(z^0)) \end{pmatrix} \right]$$

due to the formula for the Lie bracket. Since we assumed that $(l, l \circ A, ..., l \circ A^{r+1})$ is open we obtain

$$\Delta^k(z^0,...,z^{r+1}) = \sum_{j=1}^{k_D-1} \gamma_j^k(z^0,...,z^{r+1})\phi^j(z^0,...,z^r)$$

with smooth coefficients $\gamma_j^k$ by the above arguments. Differentiation with respect to $z^{r+1}$ and applying the differential to $v \in B$ yields

$$D_{z^r}\phi^k(z^0,...,z^r) \cdot v = \sum_{j=1}^{k_D-1} \left( \beta_j^k(z^0,...,z^r) \cdot v \right) \phi^j(z^0,...,z^r),$$

where $D_{z^r}\gamma_j^k(z^0,...,z^{r+1}) = \beta_j^k(z^0,...,z^r)$. Given two points $z_0^r$ and $z^r$ in $B$, we can find a smooth curve $c : \mathbb{R} \to B$ such that

$$c(0) = z_0^r \text{ and } c(1) = z^r.$$

For fixed $z^0,...,z^{r-1}$ we define $\psi^k(t) := \phi^k(z^0,...,z^{r-1},c(t))$, a smooth curve into $D(A^\infty)$, which satisfies the differential equation

$$\frac{d}{dt}\psi^k(t) = \sum_{j=1}^{k_D-1} \left( \beta_j^k(z^0,\ldots,z^{r-1},c(t))\frac{d}{dt}c(t) \right) \psi^j(t).$$

This differential equation has a unique solution, namely there exist smooth curves $A_j^k$ such that

$$\psi^k(t) = \sum_{j=1}^{k_D-1} A_j^k(t)\psi^j(0).$$

Thus for $t = 1$ there are real numbers $\alpha_j^k$ such that

$$\phi^k(z^0,...,z^r) = \sum_{j=1}^{k_D-1} \alpha_j^k \phi^k(z^0,...,z^{r-1},z_0^r).$$

By smoothness and linear independence of the fields on both sides we conclude that there are smooth matrix-valued functions $\alpha_j^k : V \times B^r \to \mathbb{R}$, with existing smooth inverse, such that

$$\phi^k(z^0,...,z^r) = \sum_{j=1}^{k_D-1} \alpha_j^k(z^0,...,z^r)\phi^k(z^0,...,z^{r-1},z_0^r).$$

Hence $\phi^k(z^0,...,z^{r-1},z_0^r)$, with frozen last variable $z_0^r$, are linearly independent and span on $V \times B^r$ the same subspace as the fields $\phi^k(z^0,...,z^r)$ do. Therefore we could have taken the $\phi^k$ with frozen last variable to span $D_{LA}$. But this is a contradiction to the minimality of $r$. Hence $r = -1$. That is, there exist constant vector fields $\Lambda_1,...,\Lambda_{k_D-1}$ around $h_0$ that span the subdistribution of $D_{LA}$ spanned by Banach map vector fields. If we are given a basis of constant vector fields around any point $h_0$, then we can extend the result to the whole domain of constant rank. This is the desired assertion. $\blacksquare$

**Remark 14.** *The time-dependent case is treated in a similar way: Generic FDRs for this invariance problem are essentially given by generic finite dimensional realizations on the extended phase space $\widetilde{D(A^\infty)} = \mathbb{R} \times D(A^\infty)$. We look at the problem on the extended phase space with operators*

$$\widetilde{\mu}(\widetilde{h}) = \begin{pmatrix} 1 \\ \mu(t,h) \end{pmatrix}$$

$$\widetilde{\sigma}^j(\widetilde{h}) = \begin{pmatrix} 0 \\ \sigma^j(t,h) \end{pmatrix}.$$

*Given an invariant manifold $\widetilde{\mathcal{M}}$ in $\widetilde{D(A^\infty)}$, the projection $pr_2 : \mathbb{R} \times D(A^\infty) \to D(A^\infty)$ restricted to $\widetilde{\mathcal{M}}$ is an immersion, since the second component of the vector fields on the extended phase space is spanned by 1. By invariance in $\widetilde{D(A^\infty)}$ we can conclude invariance in $D(A^\infty)$. Given an invariant manifold $\mathcal{M}$ in $D(A^\infty)$ we obtain a foliation by invariant manifolds in $\mathbb{R} \times \mathcal{M}$ for the homogenous problem on the extended phase space.*

*On the extended phase space we can apply the Frobenius methods. See also [3].*

**Remark 15.** *We can replace $D(A^\infty)$ in the above discussion by any Fréchet space $E$ that satisfies (H1)-(H3) (where $H$ is to be replaced by $E$). The Banach-map $\sigma : E \to E_0$ has to be chosen according to (A1)-(A3). We assume that $S(t)$ is a smooth semigroup on $E$ and (A4) has to be replaced by*

(A4') $A : E \to E$ is not a Banach map

*(see Lemma 3). We may think of $E = C^\infty(\mathbb{R}_{\geq 0}; \mathbb{R})$. Then all the above conclusions on the geometry of generic FDRs (=weak foliations) can be drawn. This illustrates that our analysis is essentially independent of the initial choice of $H$. The geometric problem even does not lead to more general solutions even on huge spaces of forward rates.*

## 5. THE SVENSSON FAMILY AS A LEAF IN A FOLIATION

A popular forward curve-fitting method is the Svensson [17] family

$$G_S(x,z) = z_1 + z_2 e^{-z_5 x} + z_3 x e^{-z_5 x} + z_4 x e^{-z_6 x}.$$

It is shown in [7] that the only non-trivial interest rate model that is consistent with the Svensson family is of the form

$$r_t = Z_t^1 g_1 + \cdots + Z_t^4 g_4, \tag{5.1}$$

where

$$g_1(x) \equiv 1, \quad g_2(x) = e^{-\alpha x}, \quad g_3(x) = x e^{-\alpha x}, \quad g_4(x) = x e^{-2\alpha x},$$

for some fixed $\alpha > 0$. Moreover,

$$Z_t^1 \equiv Z_0^1, \quad Z_t^3 = Z_0^3 e^{-\alpha t}, \quad Z_t^4 = Z_0^4 e^{-2\alpha t} \quad (Z_0^4 \geq 0)$$

and $Z^2$ satisfies

$$dZ_t^2 = \left(Z_t^3 + Z_t^4 - \alpha Z_t^2\right) dt + \sqrt{\alpha Z_t^4} dW_t. \tag{5.2}$$

Here $W$ is a real-valued Brownian motion.

We now shall find a generic local 2-dimensional realization that is of the form (5.1) whenever $r_0 = \sum_{j=1}^4 z_j g_j$ with $z_4 \geq 0$. In view of (5.2), a candidate for $\sigma$ is given, on $U := \{\ell > 0\}$, by

$$\sigma(h) = \sqrt{\alpha\ell(h)}g_2,$$

where $\ell$ is some continuous linear functional on $H$ (or even $C(\mathbb{R}_{\geq 0}, \mathbb{R})$) with $\ell(g_1) = \ell(g_2) = \ell(g_3) = 0$ and $\ell(g_4) = 1$ (notice that this is in full accordance with (4.11)). Straightforward calculations show, for $h \in U$,

$$\mu(h) = Ah + \ell(h)g_2 - \ell(h)g_2^2$$

$$[\mu, \sigma](h) = -\alpha\sqrt{\alpha\ell(h)}g_2 - \frac{\ell(\mu(h))}{2\sqrt{\alpha\ell(h)}}g_2.$$

(the clue is that $\ell \circ \sigma \equiv 0$). Hence indeed $\dim\{\mu, \sigma\}_{LA} = 2$ on $U$.

Damir Filipović, Department of Mathematics, ETH, Rämistrasse 101, CH-8092 Zürich, Switzerland. Josef Teichmann, Institute of financial and actuarial mathematics, TU Vienna, Wiedner Hauptstrasse 8-10, A-1040 Vienna, Austria

*E-mail address*: `filipo@math.ethz.ch, josef.teichmann@fam.tuwien.ac.at`